\documentclass[letterpaper,11pt,openany,oneside]{article}
\usepackage[T1]{fontenc}
\usepackage[latin1]{inputenc}
\usepackage{epsfig}
\usepackage{amsmath,amssymb}
\usepackage[active]{srcltx}
\usepackage{dsfont}
\usepackage[titles]{tocloft}
\numberwithin{equation}{section}
\newtheorem{theoreme}{Theorem}[section]

\newtheorem{definition}[theoreme]{Definition}

\newenvironment{proof}[1][Proof]{\noindent \textbf{#1.}~ }
{\hfill\rule{2mm}{2mm} \vspace{\parskip} }

\newcommand{\RR}{\ensuremath{\mathbb R}}

\newcommand{\EE}{\ensuremath{\mathbb E}}
\newcommand{\NN}{\ensuremath{\mathbb N}}

\newcommand{\I}{\ensuremath{\mathcal I}}

\newcommand{\T}{\ensuremath{\mathcal T}}
\newcommand{\J}{\ensuremath{\mathcal J}}

\newcommand{\De}{\Delta}

\def\val{\hbox{\tt val}}

\title{Uniform Value for Recursive Games with Compact Actions}
\author{Xiaoxi Li\thanks{Economics and Management School, Wuhan University, Luojia Hill, 430072 Wuhan, China. Email: xxleewhu@gmail.com.}  \ \ and 
\ Sylvain Sorin\thanks{CNRS, IMJ-PRG, UMR 7586, Sorbonne Universit\'es, UPMC Univ Paris 06, Univ Paris Diderot, Sorbonne Paris Cit\'e, Case 247, 4 Place Jussieu, 75252 Paris, France. Email: sylvain.sorin@imj-prg.fr.} } 
\date{\today}

\begin{document}
\maketitle
\abstract{Mertens, Neyman and Rosenberg \cite{Mertens_2009} used the Mertens and Neyman theorem \cite{Mertens_81} to prove the existence of uniform value for absorbing games with finite state space and compact action sets. We provide an analogous proof for another class of stochastic games, recursive games with finite state space and compact action sets. Moreover, both players have stationary $\varepsilon$-optimal strategies.} \\


\noindent \textbf{OR/MS Subject Classification}. Primary: games-group decisions/stochastic

\section{Introduction}

Zero-sum stochastic games were introduced by Shapley \cite{Shapley_53}, and the model $\Gamma=\langle K,I,J,g,q\rangle$ is as follows: $K$ is the finite state space, $I$ (resp. $J$) is the finite action set for player 1 (resp. player 2), $g:K\times I\times J\to\RR$ is the stage payoff function and $q:K\times I\times J\to \De(K)$ is the probability transition function ($\De(K)$ stands for the set of probabilities on $K$). $k_1\in K$ is the initial state. At each stage $t\geq 1$, after observing the $t$-stage history $h_t=(k_1,i_1,j_1,...,k_{t-1},i_{t-1},j_{t-1},k_t)$, player 1 chooses an action $i_t\in I$ and player 2 chooses an action $j_t\in J$. This profile $(i_t, j_t)$  induces a current stage payoff $g_t:=g(k_t,i_t,j_t)$ and transition probability $q(k_t,i_t,j_t)$ for $k_{t+1}$ the state at the next stage. 

We write $\sigma\in\Sigma$ for a behavior strategy for player 1, resp.  $\tau\in\T$ for player 2. The $\lambda$-discounted value $v_\lambda$ of the stochastic game $\Gamma$ (for $\lambda\in(0,1]$) is defined with respect to the payoff $\gamma_\lambda(k_1, \sigma,\tau)=\EE^{k_1}_{\sigma,\tau}\big[\sum_{t\geq 1} \lambda(1-\lambda)^{t-1} g(i_t, j_t, k_t)\big]$. The $n$-stage value $v_n$ (for $n\geq 1$) is defined analogously by taking the $n$-stage averaged payoff. Shapley \cite{Shapley_53} proved that
\begin {eqnarray}\label{eq:v_lambda}
 v_\lambda=\Phi(\lambda,v_\lambda), \ \text{with} 
\end {eqnarray}
\begin {eqnarray}\label{eq:Phi_lambda}
\forall f\in \RR^{|K|}, \ \  \Phi(\lambda,f)(k)=\val_{x\in\De(\I), y\in\De(\J)}\EE_{x,y}^k\big[\lambda g(i,j,k)+(1-\lambda)f(k')\big],
\end {eqnarray}
where $\val_{x\in\De(\I), y\in\De(\J)}=\max_{x\in\De(\I)}\min_{y\in\De(\J)}=\min_{y\in\De(\J)}\max_{x\in\De(\I)}$,
and that stationary optimal strategies exist for each $\lambda$, i.e. depending  at each stage $t$ only the current state $k_t$. \\ 
  
\noindent \textit{Remark}.\  When the action sets $I$ and $J$ are metric compact, one may assume $g$ and $q$ to be separately continuous on $I\times J$, and this implies the measurability for them (cf. I.1.Ex.7a in Mertens, Sorin and Zamir \cite{Mertens_2015}). The above description and results extend (cf. VII.1.a in Mertens, Sorin and Zamir \cite{Mertens_2015}): the operator $\Phi(\lambda, \cdot)$ in Equation (\ref{eq:Phi_lambda}) is well defined, $v_\lambda$ exists and satisfies Equation (\ref{eq:v_lambda}), and both players have stationary optimal strategies for each  $\lambda$. \\

We are interested in  long-run properties of $\Gamma$:  convergence of $v_\lambda$ as $\lambda$ tends to zero or  convergence of $v_n$ as $n$ tends to infinity. Moreover, in case of convergence, we ask for the existence of $\varepsilon$-optimal strategies for both players which  guarantee the limit value in all games with a sufficiently large  horizon $n$, explicitly:

\begin{definition}\label{def:epsilon-optimal} 
 \ Let $v\in \RR^{|K|}$. Play 1 can \text{guarantee} $v$ in $\Gamma$  if,  for any $\varepsilon>0$ and for every $k_1\in K$, there exist $N(\varepsilon)\in\NN$ and a behavior strategy $\sigma^*\in\Sigma$ for player 1 s.t.: \ $\forall \tau\in\T$,
\begin{eqnarray*}
&(A)& \frac{1}{n}\EE^{k_1}_{\sigma^*,\tau}\big[\sum_{t=1}^ng_t\big]\geq v(k_1)-\varepsilon,\ \forall n\geq N(\varepsilon), \cr
&(B) &  \EE^{k_1}_{\sigma^*,\tau}\big[\liminf_{n\to\infty} \frac{1}{n} \sum_{t=1}^ng_t\big]\geq v(k_1)-\varepsilon.
\end{eqnarray*}
A similar  definition holds for player 2.
\\
 $v$ is the uniform value of $\Gamma$ if both players can guarantee it. 
\end{definition}

For stochastic games with finite state space and finite action sets, Bewley and Kohlberg \cite{Bewley_76a} proved the convergence of $v_\lambda$ as $\lambda$ tends to zero (and later deduce the convergence of $v_n$ as $n$ tends to infinity), using an algebric argument. Using  the property that the function $\lambda\mapsto v_\lambda$ has bounded variation, that follows from \cite{Bewley_76a}, Mertens and Neyman \cite{Mertens_81} proved the existence of the uniform value. Actually,  Mertens and Neyman's  theorem \cite{Mertens_81} is even applicable to a stochastic game $\Gamma$ with compact action sets under the following form:


\begin{theoreme}[Mertens and Neyman \cite{Mertens_81}]\label{thm:MN}
Let $\lambda \to w_\lambda \in \RR^{|K|}$ be a function defined on $(0,1]$. Player 1 can guarantee $\limsup_{\lambda\to 0^+} w_\lambda$ in the stochastic game $\Gamma$ if $w_\lambda$ satisfies:\\
i) for some integrable function $\phi: (0,1]\to \RR_+$, $||w_\lambda -w_{\lambda'}||_\infty\leq \int_{\lambda}^{\lambda'} \phi(x)dx,\ \forall 0<\lambda,\lambda'<1$; \\
ii) for every $\lambda \in (0,1)$ sufficiently small, $\Phi(\lambda,w_\lambda) \geq  w_\lambda$.
\end{theoreme}

\noindent \textit{Remark}.\  In the construction of $\varepsilon$-optimal strategy in Mertens and Neyman \cite{Mertens_81}, $w_\lambda$ is taken to be $v_\lambda$, so condition $i)$ is implied by the bounded variation property of $v_\lambda$ and condition $ii)$ is implied by Equation (\ref{eq:v_lambda}). \\

An \textit{absorbing state} is such that once reached, the probability of leaving it is zero. Without loss of generality one  assumes that at  any absorbing state, the payoff is  absorbing (equal to the value of the static game to be played after absorption), as long as one  states that both players are informed of the current state.

$\Gamma=\langle K,I,J,g,q\rangle$ is an \textit{absorbing game} if all states but one are absorbing. Mertens, Neyman and Rosenberg \cite{Mertens_2009} used the above Theorem \ref{thm:MN} to prove the existence of uniform value for absorbing games with finite state space and compact action sets.

\textit{Recursive games}, introduced by Everett \cite{Everett_57}, are stochastic games where the stage payoff is always zero before reaching an absorbing state. This note proves the existence of uniform value for a recursive game with finite state space and compact action sets, using an approach analogous to Mertens, Neyman and Rosenberg \cite{Mertens_2009} for absorbing games. Moreover, due to the specific payoff structure, we show that $\varepsilon$-optimal strategies in recursive games can be taken stationary. This is not the case for a general stochastic game, in which an $\varepsilon$-optimal strategy has to be usually a function of the whole past history (cf. Blackwell and Ferguson \cite{Blackwell_68} for the "Big match" as an example). 

Everett \cite{Everett_57} proved the existence of stationary $\varepsilon$-optimal strategies  for both players which guarantee the "limiting-average value" (property $(B)$   in Definition \ref{def:epsilon-optimal}). As our proof relies on his characterization of this value (and on its existence), we describe here the result. \\

\noindent Given  $S\subseteq \RR^d$, let  $\overline{S}$ denote its closure.\\
 Let $K^0\subseteq K$ be the set of nonabsorbing states. \\
$\Phi(0, \cdot)$ refers to the operator $\Phi(\lambda, \cdot)$ with $\lambda=0$ in Equation (\ref{eq:Phi_lambda}). 
\begin{theoreme}[Everett \cite{Everett_57}]\label{thm:Everett} 
Let $\Gamma$ be a recursive game with finite state space and compact action sets. Then $\Gamma$ has a limiting-average value $v$ and both players have stationary $\varepsilon$-optimal strategy, i.e. $\forall \varepsilon>0$, there are stationary strategies $(\sigma^*,\tau^*)\in\Sigma\times \T$ s.t.:
for any $(\sigma,\tau)\in \Sigma\times \T$,
$$ \EE^{k_1}_{\sigma^*,\tau}\big[\liminf_{n\to\infty} \frac{1}{n} \sum_{t=1}^ng_t\big]\geq v(k_1)-\varepsilon \ \text{ and } \ \EE^{k_1}_{\sigma,\tau^*}\big[\limsup_{n\to\infty} \frac{1}{n} \sum_{t=1}^ng_t\big]\leq v(k_1)+\varepsilon$$
Moreover, the limiting-average value $v$ is characterized by\footnote{When working with the operator $\Phi(0,\cdot)$ or $\Phi(\lambda,\cdot)$, it is sufficient to consider those vectors $u\in\RR^{K}$ identical to the absorbing payoffs on the absorbing states $K\setminus K_0$. Whenever no confusion is caused, we identify $u$ with its projection on $\RR^{|K_0|}$.} $\{v\}=\overline{\xi^+} \cap \overline{\xi^-}$, where
\begin{eqnarray*}
\xi^+ = \Big\{u\in\RR^{|K^0|}: \ \Phi(0,u)\geq u, \text{ and } \Phi(0,u)(k)>u(k) \text{ whenever }  u(k)>0\Big\},\\ 
\xi^- = \Big\{u\in\RR^{|K^0|}:\ \Phi(0,u)\leq u,  \text{ and } \Phi(0,u)(k)<u(k) \text{ whenever }  u(k)<0\Big\}.
\end{eqnarray*}
\end{theoreme}      

Everett \cite{Everett_57}'s proof of the above result consists of the following two arguments: first, any vector $u\in \xi^+$ (\textit{resp.} $u\in\xi^-$) can be guaranteed by player 1 (\textit{resp.} player 2); second, the intersection of $\overline{\xi^+}$ and $\overline{\xi^-}$ is nonempty.    \\

\section{Main results and the proof}
We prove that $v$ (as characterized in Theorem \ref{thm:Everett}) is also the uniform value of $\Gamma$, and that players have stationary $\varepsilon$-optimal strategies. 
\begin{theoreme}\label{thm:main} A recursive game with finite state space and compact actions sets  has a uniform value. Moreover, both players can guarantee the uniform value in stationary strategies. 
\end{theoreme}

\noindent \textit{Remark.}  \ We emphasize that our definition of uniform value includes that of limiting-average value, thus our results extend Everett \cite{Everett_57} to a much stronger set-up. \\ 
 
\begin{proof} We first prove that $v$ is the uniform value of $\Gamma$ using Theorem \ref{thm:MN}. Let $u$ be any vector in $\xi^+$. An equivalent characterization for $u$ is: 
$$\Phi(0,u)\geq u \text{ and  } u(k)\leq 0 \text{ whenever } \Phi(0,u)(k)=u(k), \forall k\in K^0.$$ 
Define 
$w_\lambda=u, \ \forall \lambda\in(0, 1)$. We check that $w_\lambda$ satisfies the two conditions $i)$ {\&} $ii)$ in Theorem \ref{thm:MN}.

$i)$ It is trivial since $ u$ does not depend on $\lambda$;

$ii)$ The crucial point is that $\Phi(\lambda,u)= (1-\lambda)\Phi(0,u)$ for recursive games. Then the condition $\Phi(\lambda,w_\lambda)\geq w_\lambda$ for all $\lambda$ close to 0 is satisfied in the following two cases:
\begin{itemize} 
\item for $k\in K^0$ with $\Phi(0,u)(k)>u(k)$, we have $(1-\lambda)\Phi(0,u)(k)>u(k)$ for $\lambda$ close to $0^+$;
\item for $k\in K^0$ with $\Phi(0,u)(k) = u(k)$, we have $u(k)\leq0$, thus $$(1-\lambda)\Phi(0,u)(k) = (1-\lambda)u(k)\geq u(k), \ \forall \lambda\in(0,1).$$   
\end{itemize}

Now Theorem \ref{thm:MN} states that player 1 can uniformly guarantee any $u\in \xi^+$. A symmetric argument proves that  player 2 can uniformly guarantee any vector $u\in \xi^-$. As $\{v\}=\overline{\xi^+}\cap\overline{\xi^-}$ is nonempty by Theorem \ref{thm:Everett}, this proves that $v$ is the uniform value of $\Gamma$. \\ 
 
Next, we point out that the $\varepsilon$-optimal strategy appearing in the Mertens and Neyman theorem \cite{Mertens_81} can be taken \textit{stationary} when $\Gamma$ is a recursive game. Indeed, let $ w_{\lambda}$ be the function satisfying conditions $i)$ $\&$ $ii)$ in Theorem \ref{thm:MN}. The general construction of $\sigma$ for player 1 to guarantee $\limsup_{\lambda\to 0^+}w_\lambda$ is as follows: at each stage $t\geq 1$, 
\begin{itemize}
\item \underline{step 1}. the $t$-stage history is used to compute in an inductive way a small enough discount factor $\lambda_t$; 
\item \underline{step 2}. play some optimal strategy $x_{\lambda_t}(w_{\lambda_t}, k_t)\in\De(I)$ in the zero-sum game $\Phi(\lambda_t, w_{\lambda_t})(k_t)$, which is defined in Equation (\ref{eq:Phi_lambda}) by letting $f=w_{\lambda_t}$, $\lambda= \lambda_t$, $k=k_t$.  
\end{itemize}
Since 1) we have chosen in the first part of the proof $w_\lambda=u$ to be constant, 2) $\Phi(\lambda,u)= (1-\lambda)\Phi(0,u)$ for recursive games, this implies that $x_{\lambda_t}(w_{\lambda_t},k_t)$ can be taken to be some $x(u, k_t)$, an optimal strategy of player 1 in  the zero-sum game $\Phi(0, u)(k_t)=\val_{x\in\De(I), y\in\De(J)}\EE^{k_t}_{x, y}\big[ u(k')  \big]$. Hence,  to define $\sigma$, we do not need (as in Mertens and Neyman \cite{Mertens_81}) the whole $t$-stage history to compute $\lambda_t$, and the only necessary information is the current state $k_t$.   \end{proof} \\

\section{Concluding remarks}


\ \  \ \  1. To have a better understanding of the existence of stationary $\varepsilon$-optimal strategy in recursive games, one can compare our construction to the one in  Mertens, Neyman and Rosenberg \cite{Mertens_2009} for absorbing games. Indeed, they have also chosen $w_\lambda$ to be some constant function $u$. However, there is no such equality $\Phi(\lambda, u)=(1-\lambda)\Phi(0, u)$ for absorbing games, so the optimal strategy $x_{\lambda_t}(w_{\lambda_t},k_t)$ at each stage $t$ for $\Phi(\lambda_t, u)$  depends on $\lambda_t$, hence on the whole history.

2. Since the constructed $\varepsilon$-optimal strategy is stationary, our proof extends to \textit{recursive games with signals on actions}: there is no need to assume perfect observation of the opponent's actions, as long as both players are informed of the current state.    

3. A weaker notion corresponds to  the existence of an \textit{asymptotic value}, i.e. $\lim_{n\to\infty} v_n=\lim_{\lambda\to 0} v_\lambda=v$ and was studied in previous works. For absorbing games with finite state space and compact action sets, Rosenberg and Sorin \cite{Rosenberg_2001} characterized and proved such  existence, via the so-called \textit{operator approach}. Using the same approach, Sorin \cite{Sorin_2003b} provided the corresponding result for recursive games with finite state space and compact action sets and in particular,  showed that the asymptotic value  is equal to the unique vector characterized by $\overline{\xi^+}\cap\overline{\xi^-}$ in Everett \cite{Everett_57}. 

However, these convergence results do not extend to general stochastic games with finite state space and compact action sets: Vigeral \cite{Vigeral_2013} provided an example with no convergence of $v_\lambda$. This implies a fortiori that the existence result of a uniform value for absorbing/recursive games with compact actions does not extend to general  stochastic games.  

4.  Our proof does not extend to \textit{recursive games with infinite state space}. In  fact that  we use the nonemptiness of  $\overline{\xi^+}\cap \overline{\xi^-}$, which is obtained in Everett \cite{Everett_57} by an inductive proof on the (finite) number of states.\\
  Li and Venel \cite{Li_2016} provided a sufficient condition for recursive games with infinite state space to have a uniform value, that is, the family of $n$-stage values $\{v_n\}$ being totally bounded for the uniform norm.  They presented also an example to have no uniform value when this condition is not satisfied. \\ 
 
\textbf{Acknowledgement}  \ The authors are very grateful to Guillaume Vigeral for helpful comments. Part of Xiaoxi Li's work was done when he was an ATER fellow at THEMA of Universit\'e Cergy-Pontoise during the academic year 2015-2016.


\end{document}